\documentclass[10pt,a4paper]{article}
\usepackage{graphicx}
\usepackage{latexsym}
\usepackage{amssymb,amsmath,epsfig}

\newtheorem{thm}{Teorem}[section]

\newtheorem{prop}[thm]{Proposition}
\newtheorem{cor}[thm]{Corollary}
\newtheorem{defn}[thm]{Definition}
\newtheorem{ex}[thm]{Example}

\def\bm#1{\mbox{\boldmath $#1$}}
\def\proof{{\parindent 0pt \textbf{Proof: }}}

\def\qed{\, \blacksquare}

\def\I{\mathop{\rm \raisebox{0.1pt} I}}

\def\cl{\mathop{\rm cl}}

\def\atan2{\mathop{\rm atan2}}

\def\FORALL{\mbox{ for all }}

\def\thmref#1{Theorem \ref{#1}}

\def\propref#1{Proposition \ref{#1}}
\def\corref#1{Corollary \ref{#1}}

\def\exref#1{Example \ref{#1}}

\def\figref#1{Figure \ref{#1}}

\title{Domination and multistate systems}
\author{Arne Bang Huseby}

\date{}

\begin{document}
\maketitle

\begin{abstract}
\noindent Domination theory has been studied extensively in the context of binary monotone systems, where the structure function is a sum of products of the component state variables, and with coefficients given by the signed domination function. Using e.g., matroid theory, many useful properties of the signed domination function has been derived. In this paper we show how some of these results can be extended to multistate systems. In particular, we show how the signed domination function can be extended to such systems. Using M\"{o}bius inversion we show how the signed domination function can be expressed in terms of a multistate structure function. Moreover, using this expression we show how calculating the signed domination function of a multistate system can be reduced to calculating the signed domination function of an associated binary system. This way many results from binary theory can easily be extended to multistate theory.

\bigskip

\end{abstract}

\section{Introduction}
\label{sec:intro}
Binary monotone systems were introduced as a model representing how the state of a system depends on the state of the components through the structure function. For an extensive introduction to binary monotone systems and reliability calculations, see \cite{BARLOW81}. A more combinatorial approach can be found in \cite{COLBOURN87}. Even within the seemingly simple framework of binary monotone systems, reliability calculations are known to be computationally challenging. A standard method for calculating reliability is to identify the minimal path (or cut) sets, and then use the inclusion-exclusion method. By expanding the resulting expression, the structure function can be represented as a sum of products of the component state variables, and with coefficients given by the so-called signed domination function. For certain classes of systems, however, these coefficients follow certain patterns, and by utilising these, it is possible to compute the coefficients directly. This approach was used to analyse directed network systems where the coefficients are always either $+1$, $-1$ and zero. In fact for such systems the coefficients can be derived using a simple rule based on the underlying graph. For details, see \cite{SATYA78} and \cite{SATYA82}. The theory of domination plays a surprising role in other calculation methods as well, such as the factoring algorithm. For undirected network systems an optimal factoring strategy can be derived based on theoretical results related to the domination invariant. See \cite{SATYA83}. These results have later been generalised using matroid theory. See \cite{HUSEBY84}, \cite{HUSEBY89}, \cite{HUSEBY2001} and \cite{HUSEBY2011}.

An obvious limitation of binary monotone systems is that the system and the components are always described just as \emph{functioning} or \emph{failed}. In practical applications a more detailed description of the components and the system is often needed. This eventually lead to the development of the theory of multistate systems. Early attempts to this include \cite{BARLOW78}, \cite{ELNEWEIHI78}, \cite{ROSS79}, \cite{BLOCK82} and \cite{NATVIG1982}. For other early work on this, see also \cite{FUNNEMARK1985} and \cite{NATVIG1986}. Currently, the theory of multistate systems is very well established as an important branch of reliability theory, and comprehensive books on this subject include \cite{LISNIANSKI2003} and \cite{NATVIG2011}. With the recent development in machine learning methods, it has been shown that such methods can be used to approximate multistate structure functions. See \cite{ROCCO05}. This is indeed a very promising idea which surely will be developed further in the future.

Recently, an algebraic approach to multistate systems has emerged. In \cite{GIGLIO2004} a connection between multistate structure functions and monomial ideals is established. Based on this the authors show that the reliability of a multistate system can be obtained by computing the Hilbert series of the monomial ideal. An algorithm for computing this can be found in \cite{BIGATTI91}. More recent papers in this area include \cite{CABEZON2010}, \cite{CABEZON2011}, \cite{CABEZON2012}, \cite{CABEZON2014} and \cite{CABEZON2015}.

The present paper is very much inspired by the algebraic approach. However, the main goal is to show how the algebraic properties of multistate structure functions are connected to the signed domination function of binary monotone systems. This allows us to utilise the results from binary theory in the more general multistate framework. The results provided in the paper generalise many corresponding results for binary monotone systems. This also lays the groundwork for new efficient algorithms for multistate systems.

\section{Basic definitions}
\label{sec:basicDefinitions}

A \emph{multistate system} is an ordered pair $(C, \phi)$, where $C = \{1, \ldots, n\}$ is the component set, and $\phi$ is the structure function. Each component can be in a finite number of states. The set of possible states for component $i$ is denoted by $\mathcal{S}_i$, where $\mathcal{S}_i = \{0, 1, \ldots, m_i\}$, $i = 1, \ldots, n$. The component states are ranked in the natural way such that $0$ is the \emph{worst state} where the component is completely failed, while $m_i$ is the \emph{best state}, $i = 1, \ldots, n$. In order to simplify the notation we let $\mathfrak{C} = \mathcal{S}_1 \times \cdots \times \mathcal{S}_n$ denote the component state space. 

We let $x_i$ denote the state of the $i$th component, $i = 1, \ldots, n$, and let $\bm{x} = (x_1, \ldots, x_n)$ denote the \emph{component state vector}. Thus, $x_i \in \mathcal{S}_i$, $i = 1, \ldots, n$, while $\bm{x} \in \mathfrak{C}$.

Similarly, the system can be in a finite number of states, and we denote the set of possible states for the system by $\mathcal{S} = \{0, 1, \ldots, M\}$. The system states are also ranked in the natural way such that $0$ is the \emph{worst state} where the system is completely failed, and $m$ is the \emph{best state}. 

The structure function expresses the system state as a function of the component states:
\begin{align*}
\phi = \phi(\bm{x}) :  \mathfrak{C} \rightarrow \mathcal{S}
\end{align*}
For a given multistate system $(C, \phi)$ we also introduce the binary functions:
\begin{align*}
\phi_k = \phi_k(\bm{x}) = \I(\phi(\bm{x}) \geq k), \quad k = 1, \ldots, M.
\end{align*}
where $\I(\cdot)$ denotes the indicator function. We refer to $\phi_1, \ldots, \phi_M$ as the \emph{level structure functions} of the multistate system $(C, \phi)$. In particular, $\phi_k$ is called the $k$-level structure function of $(C, \phi)$, $k = 1, \ldots, M$.

A \emph{multistate monotone system} (MMS) is a multistate system $(C, \phi)$ where $\phi$ is \emph{non-decreasing} in each argument. It is easy to see that if $(C, \phi)$ is a multistate monotone system, then the level structure functions are non-decreasing in each argument as well. For convenience, the ordered pairs $(C, \phi_1), \ldots, (C, \phi_M)$ will also be referred to as multistate monotone systems, even though the structure functions only have values in $\{0,1\}$.

\section{The signed domination function}
\label{sec:signedDomination}

If $\bm{x} = (x_1, \ldots, x_n)$ and $\bm{y} = (y_1, \ldots, y_n)$, we say that $\bm{y} \leq \bm{x}$ if $y_i \leq x_i$ for all $i$. Furthermore, we say that $\bm{y} < \bm{x}$ if $\bm{y} \leq \bm{x}$ and $y_i < x_i$ for at least one $i$. If $\bm{x} \neq \bm{y}$ are such that neither $\bm{y} \leq \bm{x}$ nor $\bm{x} \leq \bm{y}$, we say that $\bm{x}$ and $\bm{y}$ are \emph{incomparable}.

If $x$ and $y$ are two numbers we denote the maximum value of $x$ and $y$ by $x \vee y$. If $\bm{x} = (x_1, \ldots, x_n)$ and $\bm{y} = (y_1, \ldots, y_n)$ are two vectors, we let:
\begin{align*}
\bm{x} \vee \bm{y} = (x_1 \vee y_1, \ldots, x_n \vee y_n)
\end{align*}

If $\mathfrak{X} = \{\bm{x}_1, \ldots, \bm{x}_s\}$ is a set of incomparable vectors, we define $\cl(\mathfrak{X})$ as the minimal set of vectors containing $\bm{x}_1, \ldots, \bm{x}_s$ which is closed under the $\vee$-operator. That is, if $\bm{x}, \bm{y} \in \cl(\mathfrak{X})$, then $\bm{x} \vee \bm{y} \in \cl(\mathfrak{X})$ as well. The set $\cl(\mathfrak{X})$ is a \emph{partially ordered set} (poset) with respect to the vector ordering defined above, where $\bm{x}_1, \ldots, \bm{x}_s$ are the minimal elements, and $\bm{x}_1 \vee \cdots \vee \bm{x}_s$ is the maximal element. We refer to $\cl(\mathfrak{X})$ as the poset \emph{generated} by the minimal elements $\bm{x}_1, \ldots, \bm{x}_s$.

It follows that if $\bm{x} \in \cl(\mathfrak{X})$, there must exist minimal vectors $\bm{x}_{i_1}, \ldots, \bm{x}_{i_j} \in \mathfrak{X}$ such that $\bm{x} = \bm{x}_{i_1} \vee \cdots \vee \bm{x}_{i_j}$. The set $\{\bm{x}_{i_1}, \ldots, \bm{x}_{i_j}\}$ is called a \emph{formation} of $\bm{x}$. If $j$ is odd, the formation is said to be \emph{odd}, while if $j$ is even, the formation is said to be \emph{even}. Note that a vector $\bm{x} \in \cl(\mathfrak{X})$ may have more than one formation. In order to keep track of these formations, the following function is useful:
\begin{align}
\label{eq:signedDomination}
\delta(\bm{x}) &= \mbox{The number of odd formations of $\bm{x}$} \\
			   &- \mbox{The number of even formations of $\bm{x}$}, \quad \bm{x} \in \cl(\mathfrak{X}) \nonumber
\end{align}
The function $\delta$ is called the \emph{signed domination function} of the poset $\cl(\mathfrak{X})$, and satisfies the following:
\begin{prop}
\label{prop:dominationInvariant}
Let $\mathfrak{X} = \{\bm{x}_1, \ldots, \bm{x}_s\}$ and $\mathfrak{Y} = \{\bm{y}_1, \ldots, \bm{y}_s\}$ be two sets of incomparable vectors such that the posets $\cl(\mathfrak{X})$ and $\cl(\mathfrak{Y})$ are isomorphic, i.e., there exists a bijective mapping $\psi : \cl(\mathfrak{X}) \rightarrow \cl(\mathfrak{Y})$ which preserves the ordering of vectors. Moreover, let $\delta$ be the signed domination function of $\cl(\mathfrak{X})$ and $\gamma$ be the signed domination function of $\cl(\mathfrak{Y})$. Then we have:
\begin{align*}
\gamma(\psi(\bm{x})) = \delta(\bm{x}), \quad \FORALL \bm{x} \in \cl(\mathfrak{X}).
\end{align*}
\end{prop}

\proof Let $\bm{x} \in \cl(\mathfrak{X})$, and assume that $\bm{x}_{i_1}, \ldots, \bm{x}_{i_j}$ is a formation of $\bm{x}$, i.e., $\bm{x} = \bm{x}_{i_1} \vee \cdots \vee \bm{x}_{i_j}$. This implies that:
\begin{align*}
\bm{x} \geq \bm{x}_{i_r}, \quad r = 1, \ldots, j.
\end{align*}
Since the mapping $\psi$ preserves the ordering of vectors, this implies that:
\begin{align*}
\psi(\bm{x}) \geq \psi(\bm{x}_{i_r}), \quad r = 1, \ldots, j.
\end{align*}
Hence, it follows that:
\begin{align}
\label{eq:inequalityA}
\psi(\bm{x}) \geq \psi(\bm{x}_{i_1}) \vee \cdots \vee \psi(\bm{x}_{i_j}).
\end{align}
We then let $\bm{y} = \psi(\bm{x}_{i_1}) \vee \cdots \vee \psi(\bm{x}_{i_j})$. This implies that:
\begin{align*}
\bm{y} \geq \psi(\bm{x}_{i_r}), \quad r = 1, \ldots, j.
\end{align*}
Since the mapping $\psi$ preserves the ordering of vectors, then so does its inverse, $\psi^{-1}$, and this implies that:
\begin{align*}
\bm{x}' = \psi^{-1}(\bm{y}) \geq \psi^{-1}(\psi(\bm{x}_{i_r})) = \bm{x}_{i_r}, \quad r = 1, \ldots, j.
\end{align*}
Hence, it follows that:
\begin{align*}
\bm{x}' \geq \bm{x}_{i_1} \vee \cdots \vee \bm{x}_{i_j} = \bm{x}
\end{align*}
Using once again that $\psi$ preserves the ordering of vectors, this implies that:
\begin{align}
\label{eq:inequalityB}
\psi(\bm{x}) \leq \psi(\bm{x}') = \bm{y} = \psi(\bm{x}_{i_1}) \vee \cdots \vee \psi(\bm{x}_{i_j})
\end{align}
By combining \eqref{eq:inequalityA} and \eqref{eq:inequalityB}, we get that:
\begin{align}
\label{eq:equalityC}
\psi(\bm{x}) = \psi(\bm{x}_{i_1}) \vee \cdots \vee \psi(\bm{x}_{i_j}).
\end{align}
Thus, if $\bm{x}_{i_1}, \ldots, \bm{x}_{i_j}$ is a formation of $\bm{x}$, then $\psi(\bm{x}_{i_1}), \ldots, \psi(\bm{x}_{i_j})$ is a formation of $\psi(\bm{x})$. Conversely, using the same arguments, it follows that if  $\psi(\bm{x}_{i_1}), \ldots, \psi(\bm{x}_{i_j})$ is a formation of $\psi(\bm{x})$, then $\bm{x}_{i_1}, \ldots, \bm{x}_{i_j}$ is a formation of $\bm{x}$. From this it follows that the number of odd formations of $\bm{x}$ is equal to the number of odd formations of $\psi(\bm{x})$, and the number of even formations of $\bm{x}$ is equal to the number of even formations of $\psi(\bm{x})$. Hence, $\gamma(\psi(\bm{x})) = \delta(\bm{x})$ as claimed $\qed$

\medskip

By \propref{eq:signedDomination} it follows that the signed domination function is an \emph{invariant} of posets. That is, if two posets are isomorphic, then their signed domination functions are equal.

\medskip

The concepts of domination and signed domination originated in the ground-breaking papers \cite{SATYA78} and \cite{SATYA83}, and other related papers from the same period. The concept was studied further in the context of matroids and oriented matroids in the papers \cite{HUSEBY84}, \cite{HUSEBY89}, \cite{HUSEBY2001} and \cite{HUSEBY2011}.

\medskip

We say that $\bm{x} \in \mathfrak{C}$ is a $k$-level \emph{path vector} of an MMS, $(C, \phi)$, if $\phi(\bm{x}) \geq k$ (or equivalently if $\phi_k(\bm{x}) = 1$). Moreover, $\bm{x}$ is a \emph{minimal} $k$-level \emph{path vector} of $(C, \phi)$ if we also have that $\phi(\bm{y}) < k$ (or equivalently if $\phi_k(\bm{y}) = 0$) for all vectors $\bm{y} \in \mathfrak{C}$ such that $\bm{y} < \bm{x}$.

For a given MMS $(C, \phi)$ we let $\mathfrak{P}_k = \{\bm{x}_1, \ldots, \bm{x}_s\}$ denote the set of minimal $k$-level path vectors. Obviously, $\phi_k(\bm{y}) = 1$ if and only if $\bm{y} \geq \bm{x}_r$ for at least one of the minimal path vectors. Thus, the $k$-level structure function, $\phi_k$ can be written as:
\begin{align}
\label{eq:unionOfPaths}
\phi_k(\bm{y}) = \I(\bigcup_{r=1}^s \{\bm{y} \geq \bm{x}_r\}), \quad \FORALL \bm{y} \in \mathfrak{C}.
\end{align}
Moreover, we note that if $\bm{x}_{i_1}, \ldots, \bm{x}_{i_j} \in \mathfrak{P}_k$, then 
\begin{align*}
\{\bm{y} \geq \bm{x}_{i_1}\} \cap \cdots \cap \{\bm{y} \geq \bm{x}_{i_j}\}
= \{\bm{y} \geq \bm{x}_{i_1} \vee \cdots \vee \bm{x}_{i_j}\}
\end{align*}
Hence, by the inclusion-exclusion formula it follows that:
\begin{align}
\label{eq:inclusionExclusion}
\phi_k(\bm{y}) &= \I(\{\bm{y} \geq \bm{x}_1\}) + \cdots + \I(\{\bm{y} \geq \bm{x}_s\}) \\[2mm]
&- \I(\{\bm{y} \geq \bm{x}_1 \vee \bm{x}_2\}) - \cdots - \I(\{\bm{y} \geq \bm{x}_{s-1} \vee \bm{x}_s\}) \nonumber \\[2mm]
&+ \I(\{\bm{y} \geq \bm{x}_1 \vee \bm{x}_2 \vee \bm{x}_3\}) 
+ \cdots + \I(\{\bm{y} \geq \bm{x}_{s-2} \vee \bm{x}_{s-1} \vee \bm{x}_s\}) \nonumber \\[2mm]
&- \cdots  \nonumber \\[2mm]
&+ (-1)^{s-1} \I(\{\bm{y} \geq \bm{x}_1 \vee \cdots \vee \bm{x}_s\}) , \quad \mbox{for all } \bm{y} \in \mathfrak{C}. \nonumber 
\end{align}
The following result shows how the signed domination function of a poset is connected to the $k$-level structure function of an MMS:
\begin{thm}
\label{thm:signedDominationA}
Let $(C, \phi)$ be an MMS, and let $\mathfrak{P}_k = \{\bm{x}_1, \ldots, \bm{x}_s\}$ denote the set of minimal $k$-level path vectors. The $k$-level structure function, $\phi_k$ can then be written as:
\begin{align}
\label{eq:signedDominationA}
\phi_k(\bm{y}) = \sum_{\bm{x} \in \cl(\mathfrak{P}_k)} \delta_k(\bm{x}) \I(\bm{y} \geq \bm{x}), \quad \FORALL \bm{y} \in \mathfrak{C}.
\end{align}
where $\delta_k$ is the signed domination function of $\cl(\mathfrak{P}_k)$, referred to as the $k$-level signed domination function of $(C, \phi)$.
\end{thm}

\proof We start out by observing that all terms in \eqref{eq:inclusionExclusion} are of the form: 
\begin{align*}
(-1)^{r-1} \I(\bm{y} \geq \bm{x}_{i_1} \vee \cdots \vee \bm{x}_{i_r}),
\end{align*}
where the odd formations have a coefficient of $(-1)^{r-1} = +1$, while the even formations have a coefficient of $(-1)^{r-1} = -1$. For each $\bm{x} \in \cl(\mathfrak{P}_k)$ we merge all terms corresponding to formations of $\bm{x}$. The merged term has a coefficient which is equal to the number of odd formations minus the number of even formations, which by definition is equal to $\delta_k(\bm{x})$. Thus, it follows that $\phi_k(\bm{y})$ satisfies \eqref{eq:signedDominationA} $\qed$

\begin{cor}
\label{cor:signedDominationA}
Let $(C, \phi)$ be an MMS, and let $\mathfrak{P}_k = \{\bm{x}_1, \ldots, \bm{x}_s\}$ denote the set of minimal $k$-level path vectors. Then we have:
\begin{align}
\label{eq:signedDominationB}
\phi_k(\bm{y}) = \sum_{\bm{x} \leq \bm{y}} \delta_k(\bm{x}), \quad \FORALL \bm{y} \in \cl(\mathfrak{P}_k),
\end{align}
where the sum in \eqref{eq:signedDominationB} is restricted to vectors $\bm{x} \in \cl(\mathfrak{P}_k)$.
\end{cor}

\proof The result follows directly from \thmref{thm:signedDominationA} by noting the obvious fact that $\I(\bm{y} \geq \bm{x}) = 1$ for all $\bm{x} \in \cl(\mathfrak{P}_k)$ such that $\bm{x} \leq \bm{y}$, and $\I(\bm{y} \geq \bm{x}) = 0$ for all $\bm{x} \in \cl(\mathfrak{P}_k)$ such that $\bm{x} \leq \bm{y}$ does not hold $\qed$

\medskip

Note that since $\phi_k(\bm{y}) = 1$ for all $\bm{y} \in \cl(\mathfrak{P}_k)$, \eqref{eq:signedDominationB} implies that:
\begin{align}
\label{eq:signedDominationC}
\sum_{\bm{x} \leq \bm{y}} \delta_k(\bm{x}) = 1, \quad \FORALL \bm{y} \in \cl(\mathfrak{P}_k),
\end{align}
where again the sum in \eqref{eq:signedDominationC} is restricted to vectors $\bm{x} \in \cl(\mathfrak{P}_k)$.

\bigskip

Before we proceed we note that if $Y_1, \ldots, Y_n$ are random variables with state spaces $\mathcal{S}_1, \ldots, \mathcal{S}_n$ respectively, and $\bm{Y} = (Y_1, \ldots, Y_n)$, then it follows by \eqref{eq:signedDominationA} that the resulting reliability of the system $(C, \phi_k)$ is given by:
\begin{align*}
P(\phi_k(\bm{Y}) = 1) \,=\, E[\phi_k(\bm{Y})] \,=\, 
\sum_{\bm{x} \in \cl(\mathfrak{P}_k)} \delta_k(\bm{x}) P(\bm{Y} \geq \bm{x})
\end{align*}
Thus, given the signed domination function $\delta_k$ and the joint distribution of $Y_1, \ldots, Y_n$, we may calculate the system reliability.

\medskip

We also note that in the binary case, i.e., where $\mathcal{S}_1 =  \cdots = \mathcal{S}_n = \{0,1\}$, we have:
\begin{align*}
\I(\bm{y} \geq \bm{x}) \,=\, \I(\bigcap_{i=1}^n y_i \geq x_i) \,=\,  \prod_{i=1}^n y_i^{x_i},
\end{align*}
where we define $y_i^0 = 1$ for all values of $y_i$, including $y_i = 0$, $i = 1, \ldots, n$. Hence, it follows that in this case \eqref{eq:signedDominationA} may alternatively be written as:
\begin{align}
\label{eq:signedDominationHilbertA}
\phi_k(\bm{y}) = \sum_{\bm{x} \in \cl(\mathfrak{P}_k)} \delta_k(\bm{x}) \prod_{i=1}^n y_i^{x_i}, \quad \FORALL \bm{y} \in \mathfrak{C}.
\end{align}
In the context of monomial ideals and Hilbert series, it is convenient to extend \eqref{eq:signedDominationHilbertA} to multistate systems. In order to avoid confusion, we use a different notation for this generalized function:
\begin{align}
\label{eq:signedDominationHilbertB}
H_k(\bm{y}) = \sum_{\bm{x} \in \cl(\mathfrak{P}_k)} \delta_k(\bm{x}) \prod_{i=1}^n y_i^{x_i}, \quad \FORALL \bm{y} \in \mathfrak{C},
\end{align}
where we again define $y_i^0 = 1$ for all values of $y_i$. See \cite{CABEZON2012}. In the binary case we obviously have $H_k(\bm{y}) = \phi_k(\bm{y})$. In the multistate case, however, this does not hold. Still there is clearly a one-to-one correspondence between $H_k(\bm{y})$ and the structure function $\phi_k(\bm{y})$. In particular, both functions share the same signed domination function. Thus, reliability calculations may just as easy be based on $H_k(\bm{y})$ as on the structure function. By using $H_k(\bm{y})$ instead of the structure function fast algorithms for Hilbert series can be applied directly in the reliability calculations. In the present paper, however, we will focus on the connection between the structure function and the signed domination function.

\bigskip

\begin{ex}
Consider an MMS $(C, \phi)$ with the set of minimal $k$-level path vectors $\mathfrak{P}_k = \{\bm{x}_1, \bm{x}_2, \bm{x}_3, \bm{x}_4\}$ where:
\begin{align*}
\bm{x}_1 = (2,1,1,0), \,\,\, \bm{x}_2 = (1,2,0,1), \,\,\,
\bm{x}_3 = (1,0,2,1), \,\,\, \bm{x}_4 = (0,1,1,2) 
\end{align*}
This implies that:
\begin{align*}
\bm{x}_1 \vee \bm{x}_2 &= (2,2,1,1), \,\,\,
\bm{x}_1 \vee \bm{x}_3 = (2,1,2,1), \,\,\,
\bm{x}_1 \vee \bm{x}_4 = (2,1,1,2), \\
\bm{x}_2 \vee \bm{x}_3 &= (1,2,2,1), \,\,\,
\bm{x}_2 \vee \bm{x}_4 = (1,2,1,2), \,\,\,
\bm{x}_3 \vee \bm{x}_4 = (1,1,2,2).
\end{align*}
Furthermore, we have:
\begin{align*}
\bm{x}_1 \vee \bm{x}_2 \vee \bm{x}_3 &= (2,2,2,1), \,\,\,
\bm{x}_1 \vee \bm{x}_2 \vee \bm{x}_4 = (2,2,1,2), \\
\bm{x}_1 \vee \bm{x}_3 \vee \bm{x}_4 &= (2,1,2,2), \,\,\,
\bm{x}_2 \vee \bm{x}_3 \vee \bm{x}_4 = (1,2,2,2),
\end{align*}
and finally, we have:
\begin{align*}
\bm{x}_1 \vee \bm{x}_2 \vee \bm{x}_3 \vee \bm{x}_4 = (2,2,2,2)
\end{align*}
Since all these $15$ vectors are distinct, it follows that each element in $\cl(\mathfrak{P}_k)$ has exactly \emph{one} formation. Thus, we get:
\begin{align*}
\delta_k(\bm{x}_{i_1}) &= +1, \quad 1 \leq i_1 \leq 4 
\\[2mm]
\delta_k(\bm{x}_{i_1} \vee \bm{x}_{i_2}) &= -1, 
\quad 1 \leq i_1 < i_2 \leq 4
\\[2mm]
\delta_k(\bm{x}_{i_1} \vee \bm{x}_{i_2} \vee \bm{x}_{i_3}) &= +1,
\quad 1 \leq i_1 < i_2 < i_3 \leq 4
\\[2mm]
\delta_k(\bm{x}_1 \vee \bm{x}_2 \vee \bm{x}_3 \vee \bm{x}_4) &= -1
\end{align*}
$\qed$
\end{ex}

\medskip

So far we have defined the signed domination function, $\delta_k$, for all $\bm{x} \in \cl(\mathfrak{P}_k)$. However, it is convenient to extend the definition of the function to the larger poset $\mathfrak{C} = \mathcal{S}_1 \times \cdots \times \mathcal{S}_n$. If $\bm{x} \notin \cl(\mathfrak{P}_k)$, obviously $\bm{x}$ does not have \emph{any formations}. Thus, by extension, \eqref{eq:signedDomination} implies that $\delta_k(\bm{x}) = 0$ for all $\bm{x} \notin \cl(\mathfrak{P}_k)$. Extending $\delta_k$ in this way ensures that we still have:
\begin{align}
\label{eq:fullDomination}
\phi_k(\bm{y}) = \sum_{\bm{x} \leq \bm{y}} \delta_k(\bm{x}), \quad \FORALL \bm{y} \in  \mathfrak{C},
\end{align}
but where the sum now ranges over \emph{all} $\bm{x} \in \mathfrak{C}$ such that $\bm{x} \leq \bm{y}$. 

\medskip
Of special interest is the signed domination of the maximal element of $\mathfrak{C}$, i.e.,  $\bm{m} = (m_1, \ldots, m_n)$. This is referred to as the \emph{signed domination} of the corresponding structure function $\phi_k$, and is denoted by $d(\phi_k)$. Thus, we have:
\begin{align}
d(\phi_k) = \delta_k(\bm{m}).
\end{align}
Note that it may happen that $\bm{m} \notin \cl(\mathfrak{P}_k)$ in which case $d(\phi_k) = 0$. This property is related to the concepts of \emph{relevance} and \emph{coherence}, well-known from the theory of binary monotone systems. In the multistate case these concepts become a bit more elaborate. Let $i \in C$, and let $r \in \{1, \ldots, m_i\}$. We say that component $i$ is \emph{$r$-relevant} with respect to $\phi_k$ if there exists some minimal path vector vector $\bm{x} \in \mathfrak{P}_k$ such that $x_i = r$. In particular, if component $i$ is $m_i$-relevant with respect to $\phi_k$, component $i$ is said to be \emph{strongly relevant} with respect to $\phi_k$. If $i$ is \emph{not} $r$-relevant for any  $r \in \{1, \ldots, m_i\}$, component $i$ is said to be \emph{irrelevant} with respect to $\phi_k$. 

The MMS $(C, \phi_k)$ is said to be \emph{strongly coherent} if all the components are strongly relevant with respect to $\phi_k$. Based on the notion of coherency, the following result is immediate:
\begin{prop}
\label{prop:coherenceDomination}
If $(C, \phi_k)$ is \emph{not} strongly coherent, then $d(\phi_k) = 0$.
\end{prop}

We observe that component relevance, as defined above, always relates to a \emph{specific} level $k$. A component may be irrelevant with respect to one level, but strongly relevant with respect to another. Note also that if component $i$ is $r$-relevant with respect to $\phi_k$, it may still happen that $i$ is \emph{not} $r'$-relevant for some $r' < r$. If we are focussing on a specific level $k$, and there are no irrelevant components, the system $(C, \phi_k)$ can be made strongly coherent simply by temporarily reducing the component state sets so that all the components become strongly relevant.

\medskip

In the binary case we obviously have $\bm{m} = \bm{1}$. This implies that the components are either $1$-relevant or irrelevant. In this case we may use the simplified term \emph{relevant} instead of $1$-relevant. Moreover, component $i$ is relevant if and only if there exists a minimal path vector $\bm{x}$ such that $x_i = 1$. Equivalently, a component is relevant if and only if it is contained in at least one minimal path set. A binary system is coherent if and only if all components are contained in at least one minimal path set.

\medskip

Note that \propref{prop:coherenceDomination} implies that strong coherency is a \emph{necessary} condition for non-zero the signed domination. As we shall see later, it is not a sufficient condition. Indeed we will present examples of strongly coherent MMS's with zero domination.

\section{M\"{o}bius inversion}
\label{sec:mobius}

The signed domination function can be analyzed using \emph{M\"{o}bius inversion}. To explain this we again consider an MMS $(C, \phi)$ and let $\mathfrak{P}_k$ denote the set of minimal $k$-level path vectors. For the poset $\cl(\mathfrak{P}_k)$ we introduce the \emph{M\"{o}bius function}, $\mu_k(\cdot, \cdot)$, defined recursively for all $\bm{x} \leq \bm{y} \in \cl(\mathfrak{P}_k)$ as follows:
\begin{align*}
\mu_k(\bm{x}, \bm{x}) &= 1, \quad \FORALL \bm{x} \in \cl(\mathfrak{P}_k), \\[2mm]
\mu_k(\bm{x}, \bm{y}) &= - \sum_{\bm{x} \leq \bm{u} < \bm{y}} \mu_k(\bm{x}, \bm{u}), \quad \FORALL \bm{x} < \bm{y} \in \cl(\mathfrak{P}_k)
\end{align*}
Then by \eqref{eq:signedDominationB} and the well-known M\"{o}bius inversion formula (see \cite{ROTA64}) we get that:
\begin{align}
\label{eq:signedDominationD}
\delta_k(\bm{y}) = \sum_{\bm{x} \leq \bm{y}} \phi_k(\bm{x}) \mu_k(\bm{x}, \bm{y}), \quad \FORALL \bm{y} \in \cl(\mathfrak{P}_k),
\end{align}
where again the sum in \eqref{eq:signedDominationD} is restricted to vectors $\bm{x} \in \cl(\mathfrak{P}_k)$. Since $\phi_k(\bm{x}) = 1$ for all $\bm{x} \in \cl(\mathfrak{P}_k)$, \eqref{eq:signedDominationD} can be simplified to:
\begin{align}
\label{eq:signedDominationE}
\delta_k(\bm{y}) = \sum_{\bm{x} \leq \bm{y}} \mu_k(\bm{x}, \bm{y}), \quad \FORALL \bm{y} \in \cl(\mathfrak{P}_k),
\end{align}

A difficulty with this approach, however, is that finding the \emph{M\"{o}bius function}, $\mu_k(\cdot, \cdot)$ for a given poset $\cl(\mathfrak{P}_k)$ is not straightforward in general. Many important results on signed domination functions can instead be obtained by working with the larger poset $\mathfrak{C}$ instead. For this poset it is possible to derive a simple explicit expression for the M\"{o}bius function. Since this function is not limited to the chosen level $k$, we denote this simply by $\mu(\cdot, \cdot)$. In order to derive $\mu(\cdot, \cdot)$ we use the following property of M\"{o}bius functions given in \cite{ROTA64}:

\begin{prop}
\label{prop:directProductRule}
Let $\mathcal{T}_1, \ldots, \mathcal{T}_r$ be finite posets with M\"{o}bius functions $\mu_1, \ldots, \mu_r$ respectively. Then the M\"{o}bius function $\mu$ of the poset $\mathcal{T}_1 \times \cdots \times  \mathcal{T}_r$ is the product of $\mu_1, \ldots, \mu_r$.
\end{prop}

In order to apply \propref{prop:directProductRule} we note that the poset $\mathfrak{C}$ is equal to the direct product $\mathcal{S}_1 \times \cdots \times  \mathcal{S}_n$. We denote the M\"{o}bius functions of $\mathcal{S}_1, \ldots,  \mathcal{S}_n$ by $\mu_1, \ldots, \mu_n$ respectively. It is easy to verify that if $x_i, y_i \in \mathcal{S}_i$ is such that $x_i \leq y_i$, $i = 1, \ldots, n$, then:
\begin{align*}
\mu_i(x_i, y_i) = \begin{cases}
 (-1)^{y_i - x_i} &\mbox{ if } x_i \leq y_i \leq x_i + 1 \\[2mm]
 0 \quad &\mbox{ otherwise} 
 \end{cases}
\end{align*}
Hence, it follows immediately by \propref{prop:directProductRule} that if $\bm{x}, \bm{y} \in \mathfrak{C}$ is such that $\bm{x} \leq \bm{y}$, then:
\begin{align}
\label{eq:mobius}
\mu(\bm{x}, \bm{y}) = \begin{cases}
 (-1)^{\sum_{i=1}^n (y_i - x_i)} &\mbox{ if } x_i \leq y_i \leq x_i + 1, \, i = 1, \ldots, n. \\[2mm]
 0 &\mbox{ otherwise} 
 \end{cases}
\end{align}
See \cite{BENDER75} for further details on this.

We can now use the M\"{o}bius function \eqref{eq:mobius} to derive an explicit expression for the signed domination function $\delta_k$. For this to work we must extend the definition of $\delta_k$ from the poset $\cl(\mathfrak{P}_k)$ to the larger poset $\mathfrak{C}$. From \eqref{eq:fullDomination} it follows by M\"{o}bius inversion that we have:
\begin{align}
\label{eq:fullMobiusInversion}
\delta_k(\bm{y}) = \sum_{\bm{x} \leq \bm{y}} \phi_k(\bm{x}) \mu(\bm{x}, \bm{y}), \quad \FORALL \bm{y} \in \mathfrak{C},
\end{align}
where the M\"{o}bius function $\mu(\cdot, \cdot)$ is given by \eqref{eq:mobius}. A more explicit formula is given in the following result:

\begin{thm}
\label{thm:signedDominationFormula}
Let $(C, \phi)$ be an MMS with component state space $\mathfrak{C}$ and system state space $\mathcal{S}$. For any $\bm{y} \in \mathfrak{C} \setminus \{\bm{0}\}$ we let $A(\bm{y}) = \{i: y_i > 0\}$. Furthermore, for any subset $B \subseteq A(\bm{y})$ we define $\bm{x}(B) = (x_1(B), \ldots, x_n(B))$ by:
\begin{align}
\label{eq:x_of_B}
x_i(B) = \begin{cases}
 y_i \quad &\FORALL i \in B \\[2mm]
 y_i - 1 \quad &\FORALL i \in A(\bm{y}) \setminus B \\[2mm]
 0 \quad &\mbox{ otherwise} 
 \end{cases}
\end{align}
Then for $k \in \{1, \ldots, M\}$, the signed domination function $\delta_k$ is given by:
\begin{align}
\label{eq:fullSignedDomination}
\delta_k(\bm{y}) = \sum_{B \subseteq A(\bm{y})} \phi_k(\bm{x}(B)) (-1)^{|A(\bm{y})| - |B|}
\end{align}
\end{thm}

\proof  We start out by noting that the M\"{o}bius function $\mu(\bm{x}, \bm{y})$ given by \eqref{eq:mobius} is non-zero if and only if $\bm{x}$ satisfies:
\begin{align}
\label{eq:nonzeroCondition}
y_i - 1 \leq x_i \leq y_i, \quad i = 1, \ldots, n.
\end{align}
Moreover, any $\bm{x}$ satisfying \eqref{eq:nonzeroCondition} is of the form $\bm{x}(B)$ where $B$ is a subset of $A(\bm{y})$. Hence, there is a one-to-one correspondence between the vectors $\bm{x} \leq \bm{y}$ such that $\mu(\bm{x}, \bm{y}) \neq 0$ and the subsets $B$ of $A(\bm{y})$. Finally, if $\bm{x} = \bm{x}(B)$ where $B \subseteq A(\bm{y})$, then:
\begin{align}
\label{eq:mobiusExponent}
\sum_{i=1}^n (y_i - x_i) = \sum_{i \in A(\bm{y})}(y_i - x_i) = |A(\bm{y}) \setminus B| = |A(\bm{y})| - |B|
\end{align}
Hence, it follows that:
\begin{align}
\label{eq:mobiusB}
\mu(\bm{x}, \bm{y}) = (-1)^{|A(\bm{y})| - |B|}, \quad \FORALL B \subseteq A(\bm{y}), \mbox{ and } \bm{x} = \bm{x}(B)
\end{align}
By inserting \eqref{eq:mobiusB} into \eqref{eq:fullMobiusInversion}, the expression \eqref{eq:fullSignedDomination} follows $\qed$

\medskip

In the \emph{binary case}, i.e., when $\mathcal{S}_i = \{0, 1\}$, $i = 1, \ldots, n$, there is a one-to-one correspondence between subsets $A$ of the component set $C$ and state vectors $\bm{x} \in \{0,1\}^n$:
\begin{align*}
A &= A(\bm{x}) = \{i : x_i = 1\}, \\[2mm]
\bm{x} &= \bm{x}(A) = (\bm{1}^{A}, \bm{0}^{C \setminus A}) 
\end{align*}
In this case the $1$-level structure function $\phi_1$ is equal to the structure function $\phi$, and consequently the signed domination function $\delta_1$ is equal to the signed domination function, $\delta$, of $\phi$. Both $\phi$ and $\delta$ can be expressed as functions of subsets of $C$, and the formula \eqref{eq:fullSignedDomination} can be written in the following simplified form:
\begin{align}
\label{eq:simplfiedSignedDomination}
\delta(A) = \sum_{B \subseteq A} \phi(B) (-1)^{|A| - |B|}, \quad 
\FORALL A \subseteq C.
\end{align}
This classical result was proved in \cite{HUSEBY84}.

\medskip

As a corollary to \thmref{thm:signedDominationFormula} we have the following result:
\begin{cor}
\label{cor:signedDominationFormula}
Let $(C, \phi)$ be an MMS with system state space $\mathcal{S}$, and where the set of states of component $i$ is $\mathcal{S}_i = \{0, 1, \ldots, m_i\}$, $i = 1, \ldots, n$. For any subset $B \subseteq C$ we define $\bm{x}(B) = (x_1(B), \ldots, x_n(B))$ by:
\begin{align*}
x_i(B) = \begin{cases}
 m_i &\FORALL i \in B \\[2mm]
 m_i - 1 &\mbox{ otherwise} 
 \end{cases}
\end{align*}
Then for $k \in \{1, \ldots, M\}$ we have:
\begin{align}
\label{eq:fullSignedDominationB}
d(\phi_k) =  \sum_{B \subseteq C} \phi_k(\bm{x}(B)) (-1)^{|C| - |B|}
\end{align}
\end{cor}

\proof The result is essentially just a special case of \thmref{thm:signedDominationFormula} where $\bm{y} = \bm{m}$, and thus $A(\bm{y}) = C$ $\qed$

\medskip

We observe that he formula for the signed domination function \eqref{eq:fullSignedDomination} given in \thmref{thm:signedDominationFormula} contains $2^{|A(\bm{y})|}$ terms. Thus, for large systems this formula is not very efficient. However, the formula turns out to be very useful for theoretical purposes. We will now use this result to prove the following \emph{recursive formula} for a signed domination function:

\begin{thm}
\label{thm:signedDominationTheorem}
Let $(C, \phi)$ be an MMS with system state space $\mathcal{S}$, and where the set of states of component $i$ is $\mathcal{S}_i = \{0, 1, \ldots, m_i\}$, $i = 1, \ldots, n$. Then for $k \in \{1, \ldots, M\}$ and $e \in C$ we have:
\begin{align}
\label{eq:signedDominationThm}
d(\phi_k) = d(\phi_k((m_e)_e, \cdot)) - d(\phi_k((m_e-1)_e, \cdot)),
\end{align}
where $\phi_k((m_e)_e, \cdot)$ denotes the structure function for the MMS obtained from $(C, \phi_k)$ by fixing the state of component $e$ to $m_e$, while $\phi_k((m_e-1)_e, \cdot)$ denotes structure function for the MMS obtained from $(C, \phi_k)$ by fixing the state of component $e$ to $m_e-1$.
\end{thm}

\proof As in \corref{cor:signedDominationFormula} we define:
\begin{align*}
x_i(B) = \begin{cases}
 m_i &\FORALL i \in B \\[2mm]
 m_i - 1 &\mbox{ otherwise} 
 \end{cases}
\end{align*}
We then use \corref{cor:signedDominationFormula} and get:
\begin{align*}
d(\phi_k) &=  \sum_{B \subseteq C} \phi_k(\bm{x}(B)) (-1)^{|C| - |B|} \\[2mm]
&= \sum_{B \subseteq C, \, e \in C} \phi_k(\bm{x}(B)) (-1)^{|C| - |B|}
+ \sum_{B \subseteq C, \, e \notin C} \phi_k(\bm{x}(B)) (-1)^{|C| - |B|} \\[2mm]
&= \sum_{B \subseteq (C\setminus e)} \phi_k(\bm{x}(B \cup e)) (-1)^{|C| - |B \cup e|} 
+ \sum_{B \subseteq (C\setminus e)} \phi_k(\bm{x}(B)) (-1)^{|C| - |B|} \\[2mm]
&= \sum_{B \subseteq (C\setminus e)} \phi_k(\bm{x}(B \cup e)) (-1)^{|C\setminus e| - |B|}
- \sum_{B \subseteq (C\setminus e)} \phi_k(\bm{x}(B)) (-1)^{|C\setminus e| - |B|} \\[2mm]
&= d(\phi_k((m_e)_e, \cdot)) - d(\phi_k((m_e-1)_e, \cdot)) \quad \qed
\end{align*}

\medskip

The formula \eqref{eq:signedDominationThm} is also referred to as a \emph{pivotal decomposition} formula for the signed domination, and the component $e \in C$, is referred to as the \emph{pivotal element}.

\medskip

The binary version of this theorem was proved in \cite{HUSEBY84}, and can be stated as a corollary of \thmref{thm:signedDominationTheorem}:
\begin{cor}
\label{cor:signedDominationTheorem}
Let $(C, \phi)$ be an binary monotone system. that is, the system state space as well as the component state spaces are $\{0,1\}$. Then for $e \in C$ we have:
\begin{align}
\label{eq:signedDominationCor}
d(\phi) = d(\phi(1_e, \cdot)) - d(\phi(0_e, \cdot)),
\end{align}
where $\phi(0_e, \cdot)$ denotes the structure function for the binary monotone system obtained from $(C, \phi)$ by fixing the state of component $e$ to $1$, while $\phi(0_e, \cdot)$ denotes structure function for the binary monotone system obtained from $(C, \phi)$ by fixing the state of component $e$ to $0$.
\end{cor}

\medskip

\begin{ex}
\label{ex:four_four_two}
Consider an MMS $(C, \phi)$ where $C = \{1,2,3,4\}$ and where:
\begin{align*}
\phi(\bm{x}) = \sum_{i=1}^4 x_i
\end{align*}
The component state spaces are $\mathcal{S}_i = \{0,1,2\}$, $i = 1,2,3,4$. From this it follows that the system state space is $\mathcal{S} = \{0, 1, \ldots, 8\}$. We then consider the $4$-level structure function:
\begin{align*}
\phi_4(\bm{x}) = \emph{I}(\sum_{i=1}^4 x_i \geq 4)
\end{align*}
Note that $(C, \phi_4)$ can be interpreted as a multistate generalization of a $k$-out-of-$n$ system, where the threshold value $4$ plays the role of $k$, but with multistate components.

The family of minimal $4$-level path vectors is $\mathfrak{P}_4 = \{\bm{x}_1, \ldots, \bm{x}_{19}\}$ where:
\begin{align*}
\bm{x}_{1} &= (2, 2, 0, 0), \quad \bm{x}_{2} = (2, 0, 2, 0), \quad 
\bm{x}_{3} = (2, 0, 0, 2), \quad \bm{x}_{4} = (0, 2, 2, 0), \\
\bm{x}_{5} &= (0, 2, 0, 2), \quad \bm{x}_{6} = (0, 0, 2, 2), \quad 
\bm{x}_{7} = (2, 1, 1, 0), \quad \bm{x}_{8} = (2, 1, 0, 1), \\
\bm{x}_{9} &= (2, 0, 1, 1), \quad \bm{x}_{10} = (1, 2, 1, 0), \quad 
\bm{x}_{11} = (1, 2, 0, 1), \quad \bm{x}_{12} = (0, 2, 1, 1), \\
\bm{x}_{13} &= (1, 1, 2, 0), \quad \bm{x}_{14} = (1, 0, 2, 1), \quad 
\bm{x}_{15} = (0, 1, 2, 1), \quad \bm{x}_{16} = (1, 1, 0, 2), \\
\bm{x}_{17} &= (1, 0, 1, 2), \quad \bm{x}_{18} = (0, 1, 1, 2), \quad 
\bm{x}_{19} = (1, 1, 1, 1)
\end{align*}

The signed domination $d(\phi_4)$ can be calculated using \eqref{eq:signedDomination}. However, this means that all possible formations must be considered. With $19$ minimal path vectors there are $2^{19}-1 = 524\,287$ such formations. Fortunately, the calculations can be simplified considerably by using the results given in this section. 

We start out by using the recursive formula \eqref{eq:signedDominationThm}. Due to the symmetry of this example the choice of pivotal element does not matter. Here we choose $e = 4$. If we fix $x_4 = 2$, the minimal $4$-level path vectors of the resulting system are:
\begin{align*}
\bm{x}_{1} &= (2, 0, 0), \quad \bm{x}_{2} = (0, 2, 0), \quad
\bm{x}_{3} = (0, 0, 2), \quad \bm{x}_{4} = (1, 1, 0), \\
\bm{x}_{5} &= (1, 0, 1), \quad \bm{x}_{6} = (0, 1, 1).
\end{align*}
If we fix $x_4 = 1$, the minimal $4$-level path vectors of the resulting system are:
\begin{align*}
\bm{x}_{1} &= (2, 1, 0),  \quad \bm{x}_{2} = (2, 0, 1),  \quad 
\bm{x}_{3} = (1, 2, 0),  \quad \bm{x}_{4} = (0, 2, 1), \\
\bm{x}_{5} &= (1, 0, 2),  \quad \bm{x}_{6} = (0, 1, 2),  \quad 
\bm{x}_{7} = (1, 1, 1)
\end{align*}
The signed dominations $d(\phi_4(2_4, \cdot))$ and $d(\phi_4(1_4, \cdot))$ can then be calculated using \eqref{eq:signedDomination} by considering a total number of $(2^{6}-1) + (2^{7}-1) = 63 + 127 = 190$ formations. Thus, we see that by using the recursive formula \eqref{eq:signedDominationThm} we have reduced the number of formations we need to consider from $524\,287$ down to just $190$ which is obviously a considerable reduction. By applying the recursive formula to the systems $(C \setminus 4, \phi_4(2_4, \cdot))$ and $(C \setminus 4, \phi_4(1_4, \cdot))$, the calculations can be further simplified.

In this case, however, there is an even simpler way of calculating $d(\phi_4)$ using \eqref{eq:fullSignedDominationB}. Since $C = \{1,2,3,4\}$ the number of terms in this sum is just $16$. Moreover, we observe that:
\begin{align*}
\bm{x}(B) \geq \bm{x}(\emptyset) = (1,1,1,1), \quad \FORALL B \subseteq C.
\end{align*}
Since $\phi_4(1,1,1,1) = 1$, it follows that $\phi_4(\bm{x}(B)) = 1$ for all $B \subseteq C$. Hence, in this case \eqref{eq:fullSignedDominationB} is reduced to:
\begin{align*}
d(\phi_k) =  \sum_{B \subseteq C} (-1)^{|C| - |B|} = 0
\end{align*}
where the last equality is a well-known combinatorial identity. Thus, we see that we have reduced the computational complexity from considering $524\,287$ formations to simply applying a well-known formula. 

Finally, note that in this example $d(\phi_4) = 0$ even though $(C, \phi_4)$ is strongly coherent $\qed$
\end{ex}

\section{Associated binary monotone systems}
\label{sec:assocBMS}

We consider the formula \eqref{eq:fullSignedDominationB} again, and note that in order to calculate $d(\phi_k)$ we only need to evaluate the function $\phi_k(\bm{x})$ for vectors $\bm{x}$ of the form $\bm{x} = \bm{x}(B)$, where:
\begin{align*}
x_i(B) = \begin{cases}
 m_i &\FORALL i \in B \\[2mm]
 m_i - 1 &\mbox{ otherwise} 
 \end{cases}
\end{align*}
By a change of variables, we can transform these vectors to binary vectors. That is, for all $B \subseteq C$ we introduce $\bm{z}(B) = (z_1(B), \ldots, z_n(B))$, where:
\begin{align*}
z_i(B) = \begin{cases}
 1 &\FORALL i \in B \\[2mm]
 0 &\mbox{ otherwise} 
 \end{cases}
\end{align*}
This implies that we have the following relation between the binary vector $\bm{z}(B)$ and the non-binary vector $\bm{x}(B)$:
\begin{align*}
\bm{x}(B) = \bm{m} - \bm{1} + \bm{z}(B), \quad \FORALL B \subseteq C.
\end{align*}
We can now associate a binary monotone system $(C, \psi_k)$ to the MMS $(C, \phi_k)$, where $\psi_k$ is a non-decreasing function defined for all binary vectors $\bm{z} \in \{0,1\}^n$ by:
\begin{align}
\label{eq:associatedBinaryStructure}
\psi_k(\bm{z}) = \phi_k(\bm{m} - \bm{1} + \bm{z})
\end{align}
We refer to $(C, \psi_k)$ as the binary monotone system \emph{associated} with the MMS $(C, \phi_k)$. By \eqref{eq:associatedBinaryStructure} we get that:
\begin{align}
\label{eq:associatedBinaryStructureB}
\psi_k(\bm{z}(B)) = \phi_k(\bm{m} - \bm{1} + \bm{z}(B)) = \phi_k(\bm{x}(B)), \quad \FORALL B \subseteq C. 
\end{align}
Inserting $\psi_k(\bm{z}(B))$ instead of $\phi_k(\bm{x}(B))$ into \eqref{eq:fullSignedDominationB}, we get:
\begin{align*}
d(\phi_k) =  \sum_{B \subseteq C} \psi_k(\bm{z}(B)) (-1)^{|C| - |B|} = d(\psi_k)
\end{align*}
The following result summarizes this discussion:
\begin{thm}
\label{thm:assocBinary}
The signed domination of a $k$-level MMS $(C, \phi_k)$ is equal to the signed domination of its associated binary monotone system $(C, \psi_k)$.
\end{thm}

More generally, similar results can be obtained for the entire signed domination function using \thmref{thm:signedDominationFormula}:
\begin{thm}
\label{thm:assocBinaryGeneral}
Let $(C, \phi)$ be an MMS with component state space $\mathfrak{C}$ and system state space $\mathcal{S}$. For any $\bm{y} \in \mathfrak{C} \setminus \{\bm{0}\}$ we let $A(\bm{y}) = \{i: y_i > 0\}$. Furthermore, for any subset $B \subseteq A(\bm{y})$ we define the subvector $\bm{z}^{A(\bm{y})}(B)$ with indices in the set $A(\bm{y})$ by:
\begin{align}
\label{eq:z_of_B}
z_i(B) = \begin{cases}
 1 \quad &\FORALL i \in B \\[2mm]
 0 \quad &\FORALL i \in A(\bm{y}) \setminus B 
 \end{cases}
\end{align}
Moreover, for $k \in \{1, \ldots, M\}$, we let $(A(\bm{y}), \psi_k^{\bm{y}})$ denote a binary monotone system where the structure function $\psi_k^{\bm{y}}$ is given by:
\begin{align*}
\psi_k^{\bm{y}}(\bm{z}^{A(\bm{y})}) = \phi_k(\bm{y} - (\bm{1}^{A(\bm{y})}, \bm{0}) + (\bm{z}^{A(\bm{y})}, \bm{0})), \quad \bm{z}^{A(\bm{y})} \in \{0,1\}^{|A(\bm{y})|}.
\end{align*}
Here $(\bm{z}^{A(\bm{y})}, \bm{0})$ denotes the binary vector where the entries with indices in the set $A(\bm{y})$ are given by the subvector $\bm{z}^{A(\bm{y})}$, while the remaining entries are zero. Similarly, $(\bm{1}^{A(\bm{y})}, \bm{0})$ denotes the binary vector where the entries with indices in the set $A(\bm{y})$ are $1$, while the remaining entries are zero.
Then we have:
\begin{align*}
\delta_k(\bm{y}) = d(\psi_k^{\bm{y}})
\end{align*}
\end{thm}

\proof We start by noting that we have the following relation between $\bm{x}(B)$ defined in \eqref{eq:x_of_B} and the subvector $\bm{z}^{A(\bm{y})}(B)$ defined in \eqref{eq:z_of_B}:
\begin{align*}
\bm{x}(B) = \bm{y} - (\bm{1}^{A(\bm{y})}, \bm{0}) + (\bm{z}^{A(\bm{y})}(B), \bm{0})), 
\quad \FORALL B \subseteq A(\bm{y}).
\end{align*}
Hence, it follows by \eqref{eq:fullSignedDomination} that we have:
\begin{align*}
\delta_k(\bm{y}) &= \sum_{B \subseteq A(\bm{y})} \phi_k(\bm{x}(B)) (-1)^{|A(\bm{y})| - |B|} \\[2mm]
&= \sum_{B \subseteq A(\bm{y})} \phi_k(\bm{y} - (\bm{1}^{A(\bm{y})}, \bm{0}) + (\bm{z}^{A(\bm{y})}(B), \bm{0}))) (-1)^{|A(\bm{y})| - |B|} \\[2mm]
&= \sum_{B \subseteq A(\bm{y})} \psi_k^{\bm{y}}(\bm{z}^{A(\bm{y})}(B)) (-1)^{|A(\bm{y})| - |B|}
= d(\psi_k^{\bm{y}})
\end{align*}
which completes the proof $\qed$

\medskip

\thmref{thm:assocBinary} and \thmref{thm:assocBinaryGeneral} are important because these results reduce the problem of computing respectively the signed domination and the signed domination function of a $k$-level MMS to the problem of computing the signed domination of associated binary monotone systems. This means that all results on signed domination of binary monotone systems can be utilized. Relevant results can be found in \cite{HUSEBY84}, \cite{HUSEBY89}, \cite{HUSEBY2001} and \cite{HUSEBY2011}, as well as in \cite{SATYA78} and \cite{SATYA83}. In the remaining part of this section we will give some examples on how \thmref{thm:assocBinary} can be used.

\subsection{Matroid systems}
\label{subsec:matroidSystems}

In \cite{HUSEBY89} and \cite{HUSEBY2001} a connection between binary monotone systems and \emph{matroids} was established. A matroid is an ordered pair $(F, \mathcal{M})$, where $F$ is a non-empty finite set, and $\mathcal{M}$ is a family of non-empty incomparable subsets of $F$, called \emph{circuits}, such that if $M_1, M_2 \in \mathcal{M}$ are two distinct circuits with a common element $e \in M_1 \cap M_2$, then there exists a third circuit $M_3 \in \mathcal{M}$ such that $M_3 \subseteq (M_1 \cup M_2) \setminus e$. See \cite{OXLEY2006} or \cite{WELSH76}.

If $(C, \phi)$ is a binary monotone system with minimal path sets $\mathcal{P}$, then $(C, \phi)$ is said to be a matroid system if there exists a matroid $(C \cup x, \mathcal{M})$ where $x \notin C$, such that:
\begin{align}
\label{eq:matroidSystem}
\mathcal{P} = \{M \setminus x : M \in \mathcal{M}, \, x \in M\}.
\end{align}
Thus, for a matroid system $(C, \phi)$ the family of minimal path sets, $\mathcal{P}$ can be derived from the matroid $(C \cup x, \mathcal{M})$ by identifying all circuits in $\mathcal{M}$ which contains the element $x$, and then delete this element from these circuits. The class of matroid systems is a large class of systems, and contains undirected network systems, $k$-out-of-$n$ systems and many other types of systems. 

If $(C, \phi)$ is a matroid system, then the matroid $(C \cup x, \mathcal{M})$ is referred to as the \emph{corresponding matroid}, and we indicate that $(C, \phi)$ can be derived from the matroid $(C \cup x, \mathcal{M})$ by writing:
\begin{align}
\label{eq:correspondingMatroid}
(C \cup x, \mathcal{M}) \rightarrow (C, \phi)
\end{align}

If $(F, \mathcal{M})$ is a matroid, a subset $B \subseteq F$ is said to be an \emph{independent set} if $B$ does not contain any circuit. The \emph{rank function} of the matroid, denoted by $\rho$, and defined for all subsets, $A \subseteq F$ is the cardinality of the largest independent subset $B \subseteq A$. That is we have:
\begin{align}
\label{eq:rankFunction}
\rho(A) = \max\{|B| : B \subseteq A \mbox{ and } M \not\subseteq B \mbox{ for all } M \in \mathcal{M}\}
\end{align}
In \cite{HUSEBY89} it was shown that $(C, \phi)$ is a matroid system if and only if the structure function $\phi$, interpreted as a function defined for all subsets $A \subseteq C$, where $A = A(\bm{x}) = \{i \in C: x_i = 1\}$, can be represented as:
\begin{align}
\label{eq:rankCharacterization}
\phi(A) = 1 + \rho(A) - \rho(A \cup x), \FORALL A \subseteq C,
\end{align}
for some matroid $(C \cup x, \mathcal{M})$ with rank function $\rho$.

The Crapo $\beta$-invariant of a matroid $(F, \mathcal{M})$ with rank function $\rho$, introduced in \cite{CRAPO67}, is given by:
\begin{align}
\label{eq:crapoBeta}
\beta(A) = \sum_{B \subseteq A} (-1)^{\rho(A) - |B|} \rho(B), \quad \FORALL B \subseteq A.
\end{align}
In particular the Crapo $\beta$-number for the matroid $(F, \mathcal{M})$, denoted $b(\mathcal{M})$ is defined as:
\begin{align}
b(\mathcal{M}) = \beta(F)
\end{align}
In \cite{CRAPO67} it was shown that $\beta(A) \geq 0$ for all $A \subseteq F$.

In \cite{HUSEBY89} it was shown that if $(C \cup x, \mathcal{M}) \rightarrow (C, \phi)$ then the signed domination function of $(C, \phi)$, denoted $\delta$, is given by:
\begin{align}
\label{eq:beta_delta}
\delta(A) = (-1)^{|A| - \rho(A \cup x)}\beta(A \cup x), \quad \FORALL A \subseteq C.
\end{align}
In particular:
\begin{align}
\label{eq:b_d}
d(\phi) = (-1)^{|C| - \rho(C \cup x)} b(\mathcal{M})
\end{align}

In \cite{SATYA83} the concept of domination was used in order to find an optimal strategy for pivotal decompositions for undirected network systems. In particular the \emph{domination invariant} of a binary monotone system $(C, \phi)$, denoted by $D(\phi)$, and defined as the absolute value of the signed domination invariant $d(\phi)$. The results on optimal pivotal decomposition strategies can be generalized to matroid systems. For details see \cite{HUSEBY89} and \cite{HUSEBY2001}. In particular, we have the following:
\begin{align}
\label{eq:Domination_partA}
D(\phi) &= b(\mathcal{M}), \\[1mm]
\label{eq:Domination_partB}
d(\phi) &= (-1)^{|C| - \rho(C \cup x)} D(\phi), \\[2mm]
\label{eq:Domination_partC}
D(\phi) &= D(\phi(1_e, \cdot)) + D(\phi(0_e, \cdot)),
\end{align}
where \eqref{eq:Domination_partC} is obtained using \eqref{eq:Domination_partB} and \corref{cor:signedDominationTheorem}.

If $(C, \phi)$ is a $k$-out-of-$n$ system where $n = |C|$ and $0 < k \leq n$, it is easy to verify that $(C, \phi)$ is a matroid system. In particular, the family of minimal path sets of $(C, \phi)$ is given by:
\begin{align*}
\mathcal{P} = \{P \subseteq C : |P| = k\},
\end{align*}
while the family of circuits of the corresponding matroid $(C \cup x, \mathcal{M})$ is given by:
\begin{align*}
\mathcal{M} = \{M \subseteq C \cup x : |M| = k+1\},
\end{align*}
The matroid $(C \cup x, \mathcal{M})$ is referred to as a \emph{uniform matroid}, and it follows by \eqref{eq:rankFunction} that $\rho(C \cup x) = k$.

Since $(C,\phi)$ is a $k$-out-of-$n$ system, it follows that $(C \setminus e, \phi(1_e, \cdot))$ is a $(k-1)$-out-of-$(n-1)$ system, while $(C \setminus e, \phi(0_e, \cdot))$ is a $k$-out-of-$(n-1)$ system. Utilizing this property in combination with \eqref{eq:Domination_partB}, \eqref{eq:Domination_partC} and induction, it is easy to verify that:
\begin{align}
\label{eq:k_out_of_n_domination}
d(\phi) = (-1)^{n-k} \binom{n-1}{k-1}, \quad 0 < k \leq n.
\end{align}

By using \thmref{thm:assocBinary} we can extend \eqref{eq:k_out_of_n_domination} to mutistate systems. One such system was analyzed in \exref{ex:four_four_two}. In the multistate case we denote such system as a $k$-out-of-$(n,m)$-system, where $k$ is a \emph{threshold value}, $n$ is the number of components in the system, and $m$ denotes a common maximal component state. More specifically, we let $(C, \phi)$ be an MMS where $C = \{1, \ldots, n\}$, $S_i = \{0, 1, \ldots, m\}$ and the structure function $\phi$ is given by:
\begin{align*}
\phi(\bm{x}) = \sum_{i=1}^n x_i
\end{align*}
Then the $k$-level MMS $(C, \phi_k)$ is said to be a $k$-out-of-$(n,m)$-system, and we have:
\begin{align*}
\phi_k(\bm{x}) = \I(\sum_{i=1}^n x_i \geq k)
\end{align*}
We then consider the associated binary monotone system $(C, \psi_k)$, where:
\begin{align*}
\psi_k(\bm{z}) &= \phi_k(\bm{m} - \bm{1} + \bm{z}) = \I(\sum_{i=1}^n (m - 1 + z_i) \geq k) \\
&= \I(\sum_{i=1}^n z_i \geq k - n(m-1))
\end{align*}
Thus, $(C, \psi_k)$ is a binary $k-n(m-1)$-out-of-$n$ system. If $k \leq n(m-1)$, it follows that $\psi_k \equiv 1$. Hence, by \eqref{eq:fullSignedDominationB} we get:
\begin{align*}
d(\phi_k) = d(\psi_k) =  \sum_{B \subseteq C} \psi_k(\bm{x}(B)) (-1)^{|C| - |B|} =  \sum_{B \subseteq C} 1 \cdot (-1)^{|C| - |B|} = 0.
\end{align*}
If $n(m-1) < k \leq nm$, it follows from \eqref{eq:k_out_of_n_domination} that:
\begin{align*}
d(\phi_k) = d(\psi_k) = (-1)^{n - (k - n(m-1))} \cdot \binom{n-1}{(k-n(m-1)) - 1} 
\end{align*}
Finally, if $k > nm$, then $\psi_k \equiv 0$. Hence, by \eqref{eq:fullSignedDominationB} we get:
\begin{align*}
d(\phi_k) = d(\psi_k) = \sum_{B \subseteq C} \psi_k(\bm{x}(B)) (-1)^{|C| - |B|} =  \sum_{B \subseteq C} 0 \cdot (-1)^{|C| - |B|} = 0.
\end{align*}
The following theorem summarizes the above discussion:
\begin{thm}
\label{thm:generalized_k_out_of_n_Domination}
Let $(C, \phi)$ be an MMS where $C = \{1, \ldots, n\}$, $S_i = \{0, 1, \ldots, m\}$ and the structure function $\phi$ is given by:
\begin{align*}
\phi(\bm{x}) = \sum_{i=1}^n x_i
\end{align*}
Then the signed domination of the $k$-level MMS $(C, \phi_k)$ is given by:
\begin{align}
\label{eq:general_k_out_of_n_domination}
d(\phi_k) = \begin{cases}
 (-1)^{n - (k - n(m-1))} \cdot \binom{n-1}{(k-n(m-1)) - 1}, &\mbox{for } n(m-1) < k \leq nm \\[1mm]
 0 &\mbox{otherwise}
 \end{cases}
\end{align}
\end{thm}

As mentioned above the MMS considered in \exref{ex:four_four_two} belongs to the same class as $(C, \phi)$ with $n = 4$ and $m = 2$. Moreover the $4$-level MMS $(C, \phi_4)$ considered in this example is a $4$-out-of-$(4,2)$-system. Since $4 \leq 4(2-1)$, it follows by the general result proved above that $d(\phi_4) = d(\psi_4) = 0$ which is the same result we arrived at in \exref{ex:four_four_two}. Note also that if $m = 1$, the general formula \eqref{eq:general_k_out_of_n_domination} reduces to \eqref{eq:k_out_of_n_domination}.

\medskip

Matroid systems also include \emph{undirected network systems}. Thus, \thmref{thm:assocBinary} can be applied to such systems as well:
\begin{ex}
\label{ex:undirectedDoubleBridgeFlow}
\begin{figure}[ht!]
  \centering
  \includegraphics[height=4.5cm]{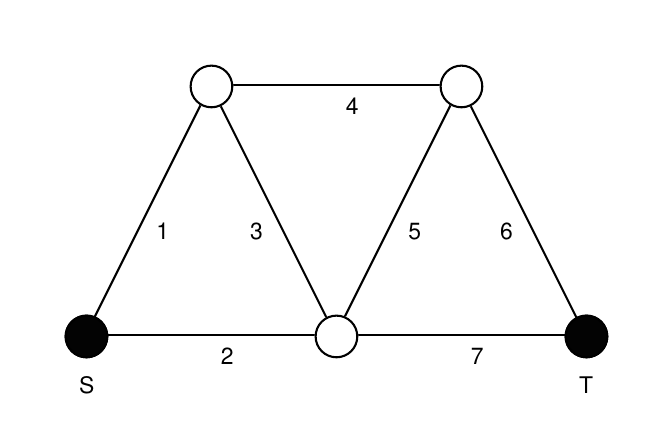}
  \caption{Undirected flow network}
  \label{fig:undirected_graph}
\end{figure}
In this example we consider the undirected 2-terminal network flow system shown in \figref{fig:undirected_graph}. This represents the MMS $(C, \phi)$, where $C = \{1, \ldots, 7\}$ is the set of edges in the network, and where $\phi(\bm{x})$ is defined as the maximum flow that can be transported through the network from the node $S$ to the node $T$. We assume that the maximum flow capacities of the edges are given by:
\begin{align*}
\bm{m} = (m_1, \ldots, m_7) = (2, 2, 1, 2, 1, 2, 2)
\end{align*}
In order to find an expression for the structure function $\phi$ we start out by finding the minimal cut sets of the network:
\begin{align*}
K_1 &= \{1,2\}, \quad K_2 = \{1,3,5,7\}, \quad K_3 = \{2,3,4\}, \\
K_4 &= \{2,3,5,6\}, \quad K_5 = \{4,5,7\}, \quad K_6 = \{6,7\}
\end{align*}
By the max-flow-min-cut theorem (see \cite{FORD62}) the structure function $\phi$ can be expressed as follows:
\begin{align*}
\phi(\bm{x}) = \min_{1 \leq r \leq 6} \sum_{i \in K_r} x_i
\end{align*}
We then consider the $3$-level MMS $(C, \phi_3)$ where $\phi_3$ is given by:
\begin{align*}
\phi_3(\bm{x}) =
\I(\sum_{i \in K_1} x_i \geq 3 \, \cap \, \cdots \, \cap \, \sum_{i \in K_6} x_i \geq 3)
\end{align*}
Hence, the structure function $\psi_3$ of the associated binary monotone system is:
\begin{align*}
\psi_k(\bm{z}) &= \phi_k(\bm{m} - \bm{1} + \bm{z}) \\[2mm]
&= \I(\sum_{i \in K_1} (m_i - 1 + z_i) \geq 3 \, \cap \, \cdots \, \cap \, \sum_{i \in K_6} (m_i - 1 + z_i) \geq 3) \\[2mm]
&= \I(\sum_{i \in K_1} z_i \geq 3 - \sum_{i \in K_1} (m_i-1)  \, \cap \, \cdots \, \cap \, \sum_{i \in K_6} z_i \geq 3 - \sum_{i \in K_1} (m_i-1)) \\[2mm]
&= \I(\sum_{i \in K_1} z_i \geq 1  \, \cap \, \cdots \, \cap \, \sum_{i \in K_6} z_i \geq 1)
\end{align*}
It is easy to verify that $\psi_k(\bm{z})$ is the structure function of the binary undirected 2-terminal network system shown in \figref{fig:undirected_graph}. In particular, $(C, \psi_k)$ is a matroid system such that the rank function $\rho$ of the corresponding matroid $(C \cup x, \mathcal{M})$ satisfies:
\begin{align*}
\rho(C \cup x) = v - 1,
\end{align*}
where $v = 5$ is the number of nodes in the graph. By using \eqref{eq:Domination_partB} it follows that:
\begin{align*}
d(\phi_3) = d(\psi_3) = (-1)^{|C| - \rho(C \cup x)} D(\psi_3) = (-1)^{7 - (5 - 1)} D(\psi_3) = - D(\psi_3)
\end{align*}
The domination $D(\psi_3)$ turns out to $3$, which is easily verified using basic properties of domination including the recursive formula \eqref{eq:Domination_partC}. Hence, it follows that $d(\phi_3) = -3$ $\qed$
\end{ex}

\subsection{Oriented matroid systems}
\label{subsec:orientedMatroidSystems}

In \cite{HUSEBY2011} the class of \emph{oriented matroid systems} was introduced. Such systems can be derived from \emph{oriented matroids} in much the same way as matroid systems are derived from matroids. For a comprehensive introduction to the theory of oriented matroids, see \cite{BJORNER99}.

Oriented matroids are based on the concept of a \emph{signed set}. A signed set is a set $M$ along with a mapping $\sigma_{M} : M \rightarrow \{+, -\}$, called the \emph{sign mapping} of the set. With a slight abuse of notation, $M$ refers both to the signed set itself as well as the underlying unsigned set of elements. The sign mapping $\sigma_{M}$ of a signed set $M$ defines a partition of $M$ into two subsets, $M^+ = \{e \in M : \sigma_{M}(e) = +\}$ and $M^- = \{e \in M : \sigma_{M}(e) = -\}$. $M^+$ and $M^-$ are referred to as the \emph{positive} and \emph{negative} elements of $M$ respectively. If $M$ is a signed set with $M^+ = \{e_1, \ldots, e_i\}$ and $M^- = \{f_1, \ldots, f_j\}$, we indicate this by writing $M$ as $\{e_1, \ldots, e_i, \bar{f}_1, \ldots, \bar{f}_j\}$. If $M = M^+$, $M$ is called a \emph{positive} set, while if $M = M^-$, $M$ is called a \emph{negative} set. An oriented matroid can now be defined as follows:
\begin{defn}
\label{defn:orientedMatroid} 
An \emph{oriented matroid} is an ordered pair $(F, \mathcal{M})$ where $F$ is a non-empty finite set, and $\mathcal{M}$ is a family of signed subsets of $F$, called \emph{signed circuits}. The signed circuits satisfy the following properties:
\begin{description}
    \item[(O1)] $\emptyset$ is not a signed circuit.
    \item[(O2)] If $M$ is a signed circuit, then so is $-M$.
    \item[(O3)] For all $M_1, M_2 \in \mathcal{M}$ such that $M_1 \subseteq M_2$, 
    		we either have $M_1 = M_2$ or $M_1 = - M_2$.
    \item[(O4)] If $M_{1}$ and $M_{2}$ are signed circuits such that $M_{1} \neq -M_{2}$, and 
    $e \in M^+_{1} \cap M^-_{2}$, then there exists a third signed circuit 
    $M_{3}$ with $M^+_{3} \subseteq (M^+_{1} \cup M^+_{2}) \setminus e$ and
    $M^-_{3} \subseteq (M^-_{1} \cup M^-_{2}) \setminus e$.
\end{description}
\end{defn}

Similar to how we defined matroid systems in \eqref{eq:matroidSystem}, we now define \emph{oriented matroid systems} as follows:
\begin{defn}
\label{defn:orientedMatroidSystem}
Let $(C \cup x, \mathcal{M})$ be an oriented matroid, and let $(C, \phi)$ be a binary monotone system with minimal path set family $\mathcal{P}$ given by:
\begin{equation}
	\label{eq:matroidSystemPaths}
    \mathcal{P} = \{(M \setminus x) : M \in \mathcal{M},\: x \in M^+ \hbox{ and } (M \setminus x)^- = \emptyset\}
\end{equation}
We then say that $(CE, \phi)$ is the oriented matroid system derived from the oriented matroid $(C \cup x, \mathcal{M})$ with respect to $x$, and write this as $(C \cup x, \mathcal{M}) \rightarrow (C, \phi)$.

If $(C \cup x, \mathcal{M}) \rightarrow (C, \phi)$, the system $(C, \phi)$ is said to be \emph{cyclic} if there exists a circuit $M \in \mathcal{M}$ such that $M \subseteq C$ and $M^- = \emptyset$. If no such circuit exists, $(C, \phi)$ is said to be \emph{acyclic}.
\end{defn}

The main results on oriented matroid systems are the following two theorems:
\begin{thm}
\label{thm:dominationAcyclicOMS}
Let $(E, \phi)$ be a coherent acyclic oriented matroid system derived from the oriented matroid $(E \cup x, \mathcal{M})$. Moreover, let $\rho$ be the rank function of $(E \cup x, \mathcal{M})$. Then $d(\phi) = (-1)^{|E| - \rho(E \cup x)}$.
\end{thm}

\begin{thm}
\label{thm:dominationCyclicOMS}
Let $(E, \phi)$ be a cyclic oriented matroid system. Then $d(\phi) = 0$.
\end{thm}

The above theorems were first proved for directed network systems in \cite{SATYA78}, while simplified proofs were obtained in \cite{HUSEBY84}. The generalization to oriented matroid systems was proved in \cite{HUSEBY2011}.

\medskip

We now illustrate how these results can be used in combination with \thmref{thm:assocBinary} by considering the following example.

\begin{ex}
\label{ex:AcyclicDirectedDoubleBridgeFlow}
\begin{figure}[ht!]
  \centering
  \includegraphics[height=4.5cm]{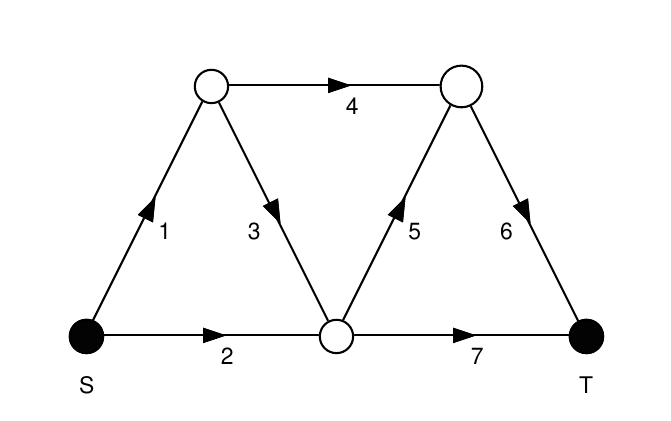}
  \caption{Acyclic directed flow network}
  \label{fig:acyclic_digraph}
\end{figure}
\figref{fig:acyclic_digraph} shows an acyclic directed 2-terminal network flow system representing the MMS $(C, \phi)$, where $C = \{1, \ldots, 7\}$ is the set of edges in the network, and where $\phi(\bm{x})$ is defined as the maximum flow that can be transported through the network from the node $S$ to the node $T$. As in \exref{ex:undirectedDoubleBridgeFlow} the maximum flow capacities of the edges are given by:
\begin{align*}
\bm{m} = (m_1, \ldots, m_7) = (2, 2, 1, 2, 1, 2, 2)
\end{align*}
The minimal cut sets of this network are:
\begin{align*}
K_1 &= \{1,2\}, \quad K_2 = \{1,5,7\}, \quad K_3 = \{2,3,4\}, \\
K_4 &= \{2,3,6\}, \quad K_5 = \{4,5,7\}, \quad K_6 = \{6,7\}
\end{align*}
Hence, by the max-flow-min-cut theorem the structure function $\phi$ is:
\begin{align*}
\phi(\bm{x}) = \min_{1 \leq r \leq 6} \sum_{i \in K_r} x_i,
\end{align*}
while the structure function of the $3$-level MMS $(C, \phi_3)$ is:
\begin{align*}
\phi_3(\bm{x}) =
\I(\sum_{i \in K_1} x_i \geq 3 \, \cap \, \cdots \, \cap \, \sum_{i \in K_6} x_i \geq 3)
\end{align*}
The calculation of the structure function $\psi_3$ of the associated binary monotone system can be done exactly in the same way as in \exref{ex:undirectedDoubleBridgeFlow}, and we eventually get:
\begin{align*}
\psi_k(\bm{z}) = \phi_k(\bm{m} - \bm{1} + \bm{z})
= \I(\sum_{i \in K_1} z_i \geq 1  \, \cap \, \cdots \, \cap \, \sum_{i \in K_6} z_i \geq 1)
\end{align*}
It is easy to verify that $\psi_k(\bm{z})$ is the structure function of the binary acyclic directed 2-terminal network system shown in \figref{fig:acyclic_digraph}. This is an acyclic oriented matroid system where the rank function $\rho$ of the oriented matroid $(C \cup x, \mathcal{M})$ satisfies:
\begin{align*}
\rho(C \cup x) = v - 1,
\end{align*}
where $v = 5$ is the number of nodes in the graph. By \thmref{thm:dominationAcyclicOMS} it follows that:
\begin{align*}
d(\phi_3) = d(\psi_3) = (-1)^{|C| - \rho(C \cup x)} = (-1)^{7-(5-1)} = -1
\end{align*}

\smallskip

\begin{figure}[ht!]
  \centering
  \includegraphics[height=4.5cm]{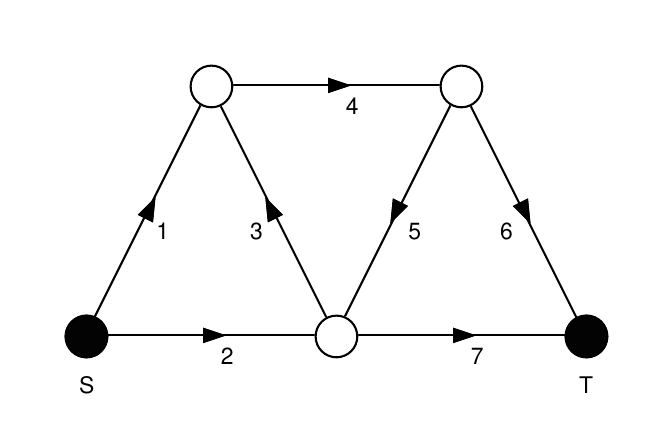}
  \caption{Cyclic directed flow network}
  \label{fig:cyclic_digraph}
\end{figure}
If we reverse the direction of the edges $3$ and $5$, we get the network shown in 
\figref{fig:cyclic_digraph}. This is a cyclic directed 2-terminal network flow system containing the directed cycle $M = M^+ = \{3,4,5\}$. The minimal cut sets then become:
\begin{align*}
K_1 &= \{1,2\}, \quad K_2 = \{1,3,7\}, \quad K_3 = \{2,4\}, \\
K_4 &= \{2,5,6\}, \quad K_5 = \{4,7\}, \quad K_6 = \{6,7\}.
\end{align*}
As before the structure function $\phi$ can be expressed as:
\begin{align*}
\phi(\bm{x}) = \min_{1 \leq r \leq 6} \sum_{i \in K_r} x_i,
\end{align*}
while the structure function of the $3$-level MMS $(C, \phi_3)$ is:
\begin{align*}
\phi_3(\bm{x}) =
\I(\sum_{i \in K_1} x_i \geq 3 \, \cap \, \cdots \, \cap \, \sum_{i \in K_6} x_i \geq 3)
\end{align*}
Finally, the structure function $\psi_3$ of the associated binary monotone system is:
\begin{align*}
\psi_k(\bm{z}) = \phi_k(\bm{m} - \bm{1} + \bm{z})
= \I(\sum_{i \in K_1} z_i \geq 1  \, \cap \, \cdots \, \cap \, \sum_{i \in K_6} z_i \geq 1)
\end{align*}
In this case, however, $\psi_k(\bm{z})$ is the structure function of the binary cyclic directed 2-terminal network system shown in \figref{fig:cyclic_digraph}. By \thmref{thm:dominationCyclicOMS} it follows that:
\begin{align*}
d(\phi_3) = d(\psi_3) = 0.
\end{align*}
Note that the signed domination is zero even though the system $(C, \phi_3)$ is strongly coherent $\qed$
\end{ex}

\section{Conclusions and further work}
\label{sec:conclusions}
We have studied the relation between $k$-level structure functions of an MMS and the corresponding \emph{signed domination function}. Formally, the signed domination function is defined through the concept of \emph{formations} of a poset of vectors. Computing the signed domination function directly from formations is very time consuming. By using \emph{M\"o{}bius inversion} the signed domination function of a $k$-level structure function can be expressed in terms of this structure function. The inversion formula shows that the signed domination of a $k$-level structure function is equal to the signed domination of its associated binary structure function. This result implies that the theory of signed domination functions for binary monotone systems can be applied to MMSs as well. Future work on this topic will further explore the use of such results.

\bibliographystyle{plain}
\bibliography{lib_file}

\end{document}